\documentclass{amsart}
\usepackage{amssymb}
\usepackage{amscd}
\usepackage{verbatim}
\usepackage{epsfig}

\begin{document}
\newcommand\Mand{\ \text{and}\ }
\newcommand\Mor{\ \text{or}\ }
\newcommand\Mfor{\ \text{for}\ }
\newcommand\Real{\mathbb{R}}
\newcommand\RR{\mathbb{R}}
\newcommand\im{\operatorname{Im}}
\newcommand\re{\operatorname{Re}}
\newcommand\sign{\operatorname{sign}}
\newcommand\sphere{\mathbb{S}}
\newcommand\BB{\mathbb{B}}
\newcommand\HH{\mathbb{H}}
\newcommand\dS{\mathrm{dS}}
\newcommand\ZZ{\mathbb{Z}}
\newcommand\codim{\operatorname{codim}}
\newcommand\Sym{\operatorname{Sym}}
\newcommand\End{\operatorname{End}}
\newcommand\Span{\operatorname{span}}
\newcommand\Ran{\operatorname{Ran}}
\newcommand\ep{\epsilon}
\newcommand\Cinf{\cC^\infty}
\newcommand\dCinf{\dot \cC^\infty}
\newcommand\CI{\cC^\infty}
\newcommand\dCI{\dot \cC^\infty}
\newcommand\Cx{\mathbb{C}}
\newcommand\Nat{\mathbb{N}}
\newcommand\dist{\cC^{-\infty}}
\newcommand\ddist{\dot \cC^{-\infty}}
\newcommand\pa{\partial}
\newcommand\Card{\mathrm{Card}}
\renewcommand\Box{{\square}}
\newcommand\Ell{\mathrm{Ell}}
\newcommand\WF{\mathrm{WF}}
\newcommand\WFh{\mathrm{WF}_\semi}
\newcommand\WFb{\mathrm{WF}_\bl}
\newcommand\WFsc{\mathrm{WF}_\scl}
\newcommand\Vf{\mathcal{V}}
\newcommand\Vb{\mathcal{V}_\bl}
\newcommand\Vz{\mathcal{V}_0}
\newcommand\Hb{H_{\bl}}
\newcommand\Hsc{H_{\scl}}
\newcommand\Ker{\mathrm{Ker}}
\newcommand\Range{\mathrm{Ran}}
\newcommand\Hom{\mathrm{Hom}}
\newcommand\Id{\mathrm{Id}}
\newcommand\sgn{\operatorname{sgn}}
\newcommand\ff{\mathrm{ff}}
\newcommand\tf{\mathrm{tf}}
\newcommand\esssupp{\operatorname{esssupp}}
\newcommand\supp{\operatorname{supp}}
\newcommand\vol{\mathrm{vol}}
\newcommand\Diff{\mathrm{Diff}}
\newcommand\Diffd{\mathrm{Diff}_{\dagger}}
\newcommand\Diffs{\mathrm{Diff}_{\sharp}}
\newcommand\Diffb{\mathrm{Diff}_\bl}
\newcommand\DiffbI{\mathrm{Diff}_{\bl,I}}
\newcommand\Diffbeven{\mathrm{Diff}_{\bl,\even}}
\newcommand\Diffz{\mathrm{Diff}_0}
\newcommand\Psih{\Psi_{\semi}}
\newcommand\Psihcl{\Psi_{\semi,\cl}}
\newcommand\Psisc{\Psi_\scl}
\newcommand\Psib{\Psi_\bl}
\newcommand\Psibc{\Psi_{\mathrm{bc}}}
\newcommand\TbC{{}^{\bl,\Cx} T}
\newcommand\Tb{{}^{\bl} T}
\newcommand\Sb{{}^{\bl} S}
\newcommand\Tsc{{}^{\scl} T}
\newcommand\Ssc{{}^{\scl} S}
\newcommand\Lambdab{{}^{\bl} \Lambda}
\newcommand\zT{{}^{0} T}
\newcommand\Tz{{}^{0} T}
\newcommand\zS{{}^{0} S}
\newcommand\dom{\mathcal{D}}
\newcommand\cA{\mathcal{A}}
\newcommand\cB{\mathcal{B}}
\newcommand\cE{\mathcal{E}}
\newcommand\cG{\mathcal{G}}
\newcommand\cH{\mathcal{H}}
\newcommand\cU{\mathcal{U}}
\newcommand\cO{\mathcal{O}}
\newcommand\cF{\mathcal{F}}
\newcommand\cM{\mathcal{M}}
\newcommand\cQ{\mathcal{Q}}
\newcommand\cR{\mathcal{R}}
\newcommand\cI{\mathcal{I}}
\newcommand\cL{\mathcal{L}}
\newcommand\cK{\mathcal{K}}
\newcommand\cC{\mathcal{C}}
\newcommand\cX{\mathcal{X}}
\newcommand\cY{\mathcal{Y}}
\newcommand\cP{\mathcal{P}}
\newcommand\cS{\mathcal{S}}
\newcommand\cZ{\mathcal{Z}}
\newcommand\cW{\mathcal{W}}
\newcommand\Ptil{\tilde P}
\newcommand\ptil{\tilde p}
\newcommand\chit{\tilde \chi}
\newcommand\yt{\tilde y}
\newcommand\zetat{\tilde \zeta}
\newcommand\xit{\tilde \xi}
\newcommand\taut{{\tilde \tau}}
\newcommand\phit{{\tilde \phi}}
\newcommand\mut{{\tilde \mu}}
\newcommand\sigmat{{\tilde \sigma}}
\newcommand\sigmah{\hat\sigma}
\newcommand\zetah{\hat\zeta}
\newcommand\etah{\hat\eta}
\newcommand\loc{\mathrm{loc}}
\newcommand\compl{\mathrm{comp}}
\newcommand\reg{\mathrm{reg}}
\newcommand\GBB{\textsf{GBB}}
\newcommand\GBBsp{\textsf{GBB}\ }
\newcommand\bl{{\mathrm b}}
\newcommand\scl{{\mathrm{sc}}}
\newcommand{\sH}{\mathsf{H}}
\newcommand{\cte}{\digamma}
\newcommand\cl{\operatorname{cl}}
\newcommand\hsf{\mathcal{S}}
\newcommand\Div{\operatorname{div}}
\newcommand\hilbert{\mathfrak{X}}
\newcommand\smooth{\mathcal{J}}
\newcommand\decay{\ell}
\newcommand\symb{j}

\newcommand\Hh{H_{\semi}}

\newcommand\bM{\bar M}
\newcommand\Xext{X_{-\delta_0}}

\newcommand\xib{{\underline{\xi}}}
\newcommand\etab{{\underline{\eta}}}
\newcommand\zetab{{\underline{\zeta}}}

\newcommand\xibh{{\underline{\hat \xi}}}
\newcommand\etabh{{\underline{\hat \eta}}}
\newcommand\zetabh{{\underline{\hat \zeta}}}

\newcommand\zn{z}
\newcommand\sigman{\sigma}
\newcommand\psit{\tilde\psi}
\newcommand\rhot{{\tilde\rho}}

\newcommand\hM{\hat M}

\newcommand\Op{\operatorname{Op}}
\newcommand\Oph{\operatorname{Op_{\semi}}}

\newcommand\innr{{\mathrm{inner}}}
\newcommand\outr{{\mathrm{outer}}}
\newcommand\full{{\mathrm{full}}}
\newcommand\semi{\hbar}

\newcommand\Feynman{\mathrm{Feynman}}
\newcommand\future{\mathrm{future}}
\newcommand\past{\mathrm{past}}

\newcommand\elliptic{\mathrm{ell}}
\newcommand\diffordgen{k}
\newcommand\difford{2}
\newcommand\diffordm{1}
\newcommand\diffordmpar{1}
\newcommand\even{\mathrm{even}}
\newcommand\dimn{n}
\newcommand\dimnpar{n}
\newcommand\dimnm{n-1}
\newcommand\dimnp{n+1}
\newcommand\dimnppar{(n+1)}
\newcommand\dimnppp{n+3}
\newcommand\dimnppppar{n+3}

\newcommand\sob{s}

\newtheorem{lemma}{Lemma}
\newtheorem{prop}[lemma]{Proposition}
\newtheorem{thm}[lemma]{Theorem}
\newtheorem{cor}[lemma]{Corollary}
\newtheorem{result}[lemma]{Result}
\newtheorem*{thm*}{Theorem}
\newtheorem*{prop*}{Proposition}
\newtheorem*{cor*}{Corollary}
\newtheorem*{conj*}{Conjecture}
\theoremstyle{remark}
\newtheorem{rem}[lemma]{Remark}
\newtheorem*{rem*}{Remark}
\theoremstyle{definition}
\newtheorem{Def}[lemma]{Definition}
\newtheorem*{Def*}{Definition}

\newcommand{\mar}[1]{{\marginpar{\sffamily{\scriptsize #1}}}}
\newcommand\av[1]{\mar{AV:#1}}

\renewcommand{\theenumi}{\roman{enumi}}
\renewcommand{\labelenumi}{(\theenumi)}

\title{On the positivity of propagator differences}
\author[Andras Vasy]{Andr\'as Vasy}
\address{Department of Mathematics, Stanford University, CA 94305-2125, USA}

\email{andras@math.stanford.edu}


\thanks{The author gratefully
  acknowledges partial support from the NSF under grant numbers
  DMS-1068742 and  DMS-1361432.}

\begin{abstract}
We discuss positivity properties of `distinguished propagators', i.e.\
distinguished {\em inverses} of operators that frequently occur in
scattering theory and wave propagation. We relate this to the work of
Duistermaat and H\"ormander on distinguished {\em parametrices} (approximate
inverses), which has played a major role in quantum field theory on
curved spacetimes recently.
\end{abstract}

\maketitle

\section{Introduction and main results}
In this short paper we discuss positivity properties of the differences of
`propagators', i.e.\ inverses of operators of the kind that frequently
occur in scattering theory and wave propagation. Concretely, we
discuss various settings in which there are function spaces,
corresponding to the `distinguished parametrices' of Duistermaat and
H\"ormander \cite{FIOII}, on which
these operators are Fredholm; in the case of actual invertibility one
has inverses and one can ask about the positivity properties of their
differences. As we recall below, Duistermaat and H\"ormander analyzed
possibilities for choices of {\em parametrices} (approximate inverses
modulo smoothing) possessing such positivity
properties; here we show that certain of the {\em actual inverses} possess
these properties, and we give a new proof of the
Duistermaat-H\"ormander theorem when our Fredholm setup is applicable.
Such a result is relevant to quantum field theory on curved
spacetimes, with work in this direction, relying on the
Duistermaat-H\"ormander framework, initiated by Radzikowski
\cite{Radzikowski:Microlocal}; see the work of Brunetti, Fredenhagen and K\"ohler
\cite{Brunetti-Fredenhagen-Kohler:Microlocal,
  Brunetti-Fredenhagen:Microlocal}, of Dappiaggi, Moretti and Pinamonti
\cite{Dappiaggi-Moretti-Pinamonti:Rigorous, Moretti:Quantum,
  Dappiaggi-Moretti-Pinamonti:Cosmological} and of G\'erard and Wrochna
\cite{Gerard-Wrochna:Construction, Gerard-Wrochna:Yang-Mills} for more
recent developments. It turns out that the positivity properties are
closely connected to the positivity of spectral measure for the
Laplacian in scattering theory.

As background, we
first recall that in elliptic settings,
or microlocally (in $T^*X\setminus o$) where a pseudodifferential
operator $P$ on a manifold $X$ is elliptic,
there are no choices to make: parametrices (as well as inverses when
one has a globally well-behaved `fully elliptic' problem and these
exist) are essentially unique; here for parametrices uniqueness is up to smoothing terms. On the other hand, if $P$ is scalar with real
principal symbol $p$ (with a homogeneous representative), or simply has
real scalar
principal symbol, then H\"ormander's theorem
\cite{Hormander:Existence} states that singularities of solutions to
$Pu=f$ propagate along bicharacteristics (integral curves of the
Hamilton vector field $H_p$) in the characteristic set $\Sigma$, in
the sense that
$\WF^{\sob}(u)\setminus\WF^{\sob-m+1}(Pu)\subset\Sigma$ is
invariant under the Hamilton flow; here $m$ is the order of
$P$. In terms of estimates, the propagation theorem states that one
can estimate $u$ in $H^\sob$ microlocally at a point $\alpha\in T^*X\setminus o$
if one has an a priori estimate for $u$ in $H^\sob$ at $\gamma_\alpha(t)$ for some
$t>0$, where $\gamma_\alpha$ is the bicharacteristic through $\alpha$,
and if one has an a priori estimate for $Pu$ in $H^{\sob-m+1}$
microlocally along $\gamma_\alpha|_{[0,t]}$; the analogous statement
for $t<0$ also holds.

Such a propagation statement is empty where $H_p$ is radial,
i.e.\ is a multiple of the radial vector field in $T^*X\setminus o$,
with the latter being the infinitesimal generator of dilations in the
fibers of $T^*X\setminus o$. However, these radial points also have
been analyzed, starting with the work of Guillemin and Schaeffer
\cite{Guillemin-Schaeffer:Fuchsian} in the case of isolated radial
points, further explored by Hassell, Melrose and Vasy
\cite{Hassell-Melrose-Vasy:Spectral,Hassell-Melrose-Vasy:Microlocal}
inspired by the work of Herbst \cite{Herbst:Spectral} and Herbst and
Skibsted \cite{Herbst-Skibsted:Absence} on a scattering problem,
by Melrose \cite{RBMSpec} for Lagrangian submanifolds of normal
sources/sinks in scattering theory,
and by Vasy \cite{Vasy-Dyatlov:Microlocal-Kerr} in a very general
situation (more general than radial points), with a more detailed analysis by Haber and Vasy in
\cite{Haber-Vasy:Radial}; see also the work of Dyatlov and Zworski
\cite{Dyatlov-Zworski:Dynamical} for their role in dynamical systems. (In a more complicated direction, in
$N$-body scattering these correspond to the propagation set of Sigal
and Soffer \cite{Sigal-Soffer:N}; see \cite{GIS:N-body} for a
discussion that is microlocal in the radial variable and see
\cite{Vasy:Bound-States} for a fully
microlocal discussion.)

In order to make the picture very clear,
consider the Hamilton flow on $S^*X=(T^*X\setminus o)/\RR^+$ rather than on $T^*X\setminus o$.
This is possible if $m=1$ since the Hamilton vector field then is
homogeneous of degree $0$ and thus can be thought of as a vector field
on $S^*X$. For general $m$ one can reduce to this case by multiplying
by a positive elliptic factor; the choice of the elliptic factor changes the
Hamilton vector field but within $\Sigma$ only by a positive factor;
in particular the bicharacteristics only get reparameterized. Thus a
radial point is a critical point for the Hamilton vector field on
$S^*X$ (i.e.\ where the vector field vanishes); in the cases discussed
here it is a non-degenerate source or sink.

In fact, it is better to
think of $S^*X$ as `fiber infinity' $\pa\overline{T^*}X$ on the fiber
compactification of $T^*X$. Here recall that if $V$ is a $k$-dimensional vector space,
it has a natural compactification $\overline{V}$ obtained by gluing a
sphere, namely $(V\setminus 0)/\RR^+$ to infinity. Explicitly this can be done
e.g.\ by putting a (positive definite) inner product on
$V$, so $V\setminus 0$ is identified with
$\RR^+_r\times\sphere^{k-1}$, with $r$ the distance from $0$, and using
`reciprocal polar coordinates'
$(\rho,\omega)\in(0,\infty)\times\sphere^{k-1}$, $\rho=r^{-1}$, to glue
in the sphere at $\rho=0$, so that the resulting manifold is covered
with the two (generalized) coordinate charts $V$ and
$[0,\infty)_\rho\times\sphere^{k-1}$ with overlap $V\setminus 0$,
resp.\ $(0,\infty)_\rho\times\sphere^{k-1}$, identified as above. This process gives a smooth
structure independent of choices, and correspondingly it can be
applied to compactify the fibers of $T^*X$. For standard microlocal
analysis the relevant location is fiber infinity, so one may instead
simply work with $S^*X\times[0,\ep)_\rho$, if one so desires, with
the choice of a homogeneous degree $-1$ function $\rho$ on
$T^*X\setminus o$ giving the identification.

The advantage for this point of view is that the
Hamilton vector field in fact induces a vector field
$\sH_p=\rho^{m-1}H_p$ on
$\overline{T^*}X$, tangent to $\pa \overline{T^*}X$, whose
linearization at radial points in $\pa\overline{T^*}X$ is well
defined. This includes the normal to the fiber boundary behavior,
i.e.\ that on homogeneous degree $-1$ functions on $T^*X\setminus o$,
via components $\rho\pa_\rho$ of the vector field; this disappears in the quotient picture. We are then interested in
critical points that are sources/sinks within $\Sigma$ even in this extended sense, so
$H_p\rho=\rho^{-m+2}\beta_0$, where $\rho$ is a boundary defining
function, e.g.\ a positive homogeneous degree $-1$ function on
$T^*X\setminus o$ near $\pa\overline{T^*}X$, and where $\beta_0>0$ at
sources, $\beta_0<0$ at sinks.
Such behavior is automatic for Lagrangian submanifolds of radial
points (these are the maximal dimensional sets of non-degenerate
radial points). The typical basic result is that there is a
threshold regularity $\sob_0$ such that for $\sob<\sob_0$ one
has a propagation of singularities type result: if a punctured
neighborhood $U\setminus\Lambda$ of a source/sink type radial set $\Lambda$ is disjoint from
$\WF^\sob(u)$ and the corresponding neighborhood $U$ is disjoint from
$\WF^{\sob-m+1}(Pu)$, then $\Lambda\cap\WF^\sob(u)=\emptyset$,
i.e.\ one can propagate estimates into $\Lambda$, while if
$\sob>\sob_1>\sob_0$, and
$\WF^{\sob_1}(u)\cap\Lambda=\emptyset$ then one can gets {\em `for
  free'} $H^\sob$ regularity at $\Lambda$, i.e.\
$\WF^\sob(u)\cap\Lambda=\emptyset$.

Here we emphasize that all of the results below hold in the more
general setting discussed in
\cite[Section~2.2]{Vasy-Dyatlov:Microlocal-Kerr}, where $\Lambda$ are
`normal sources/sinks', but need not consist of actual radial points,
i.e.\ there may be a non-trivial Hamilton flow within $\Lambda$ ---
this is the case for instance in problems related to Kerr-de Sitter
spaces. Furthermore, the setup is also stable under general
pseudodifferential (small!) perturbations of order $m$ (with real
principal symbol), even though the dynamics can change under these;
this is due to the stability of the estimates (and the corresponding
stability of the normal dynamics
in a generalized sense)
see \cite[Section~2.7]{Vasy-Dyatlov:Microlocal-Kerr}.

Now, the estimates given by the propagation theorem let one estimate
$u$ somewhere in terms of $Pu$ {\em provided one has an estimate for
  $u$ somewhere else}. But where can such an estimate come from? A
typical situation for hyperbolic equations is Cauchy data, which is
somewhat awkward from the microlocal analysis perspective and indeed
is very ill-suited to Feynman type propagators. A more natural place
is from radial sets: if one is in a sufficiently regular (above the
threshold) Sobolev space, one gets regularity for free there in terms
of a weaker (but stronger than the threshold) Sobolev norm. (This
weaker norm is relatively compact in the settings of interest, and
thus is irrelevant for Fredholm theory.) This can
then be propagated along bicharacteristics, and indeed can be
propagated into another radial set provided that we use Sobolev spaces
which are weaker than the threshold regularity there. This typically
requires the use of variable order Sobolev spaces, but as the
propagation of singularities still applies for these, provided the
Sobolev order is monotone decreasing in the direction in which we
propagate our estimates (see \cite[Appendix]{Baskin-Vasy-Wunsch:Radiation}), this is not a problem. Note that in order to
obtain Fredholm estimates eventually we need analogous estimates for
the adjoint (relative to $L^2$) $P^*$ on dual (relative to $L^2$)
spaces; since the dual of above, resp.\ below threshold regularity is
regularity below, resp.\ above threshold regularity, for the adjoint
one will need to propagate estimates in the {\em opposite direction}.
Notice that {\em within each connected component} one has to have the same
direction of propagation relative to the Hamilton flow, but of course
one can make different choices in different connected components. This
general framework was introduced by the author in
\cite{Vasy-Dyatlov:Microlocal-Kerr}, further developed with Baskin and
Wunsch in \cite{Baskin-Vasy-Wunsch:Radiation}, with Hintz in
\cite{Hintz-Vasy:Semilinear} and with Gell-Redman and Haber in \cite{Gell-Redman-Haber-Vasy:Feynman}.

Returning to the main theme of the paper, we recall that in their
influential paper \cite{FIOII} Duistermaat and H\"ormander used the
Fourier integral operators they just developed to
construct distinguished parametrices for real principal type
equations: for each component of the characteristic set, one chooses
the direction in which estimates, or equivalently singularities of
forcing (i.e.\ of $f$ for $u$ being the parametrix applied to $f$) propagate along the Hamilton flow in
the sense discussed above. Here the direction is most
conveniently measured relative to the Hamilton flow in the
characteristic set. Thus, with $k$ components of the characteristic
set, there are $2^k$ distinguished parametrices. Notice that there are
two special choices for the distinguished parametrices: the one
propagating estimates forward everywhere along the $H_p$-flow, and the
one propagating estimates backward everywhere along the $H_p$-flow;
these are the Feynman and anti-Feynman parametrices (defined up to
smoothing operators). Duistermaat and H\"ormander
showed that, if
the operator $P$ is formally self-adjoint, one
can choose these parametrices (which are defined modulo smoothing
operators a priori) so that they are all formally
skew-adjoint, and further such that $\imath$ times the difference
between any of these parametrices and the Feynman, i.e.\ the $H_p$-forward, parametrix is
positive. They also stated that they do not
see a way of fixing the smoothing ambiguity, though the paper suggests
that this would be important in view of the relationship to quantum
field theory, as suggested to the authors by Wightmann.

The purpose of this paper is to show how, under a natural additional
assumption on the global dynamics, the ambiguity can be fixed for all propagators, and
exact positivity can be shown for the extreme difference of
propagators. A byproduct is a simple proof of the positivity for a
suitable choice of distinguished parametrices (not just the extreme
difference), giving a different proof of the Duistermaat-H\"ormander
result. However, one cannot expect in general that the differences other than the
extreme difference are actually positive; thus, if positivity is
desired, the {\em only} natural choice is that of the Feynman propagators.

\begin{figure}[ht]
\includegraphics[width=60mm]{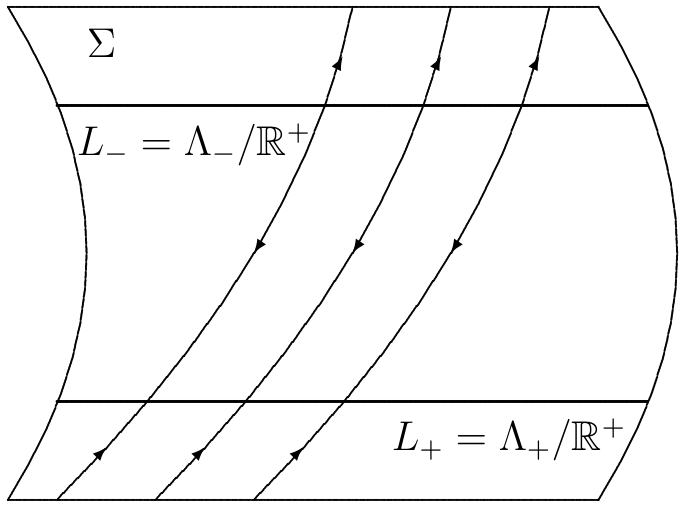}
\caption{The characteristic set $\Sigma$ (here connected) and the Hamilton dynamics for
  a problem satisfying our assumptions. Here $\Sigma$ is a torus, with
  the left and right, as well as the top and bottom, edges
  identified. An example is the
  multiplication operator by a real valued function on a compact
  manifold with non-degenerate zeros. This is closely related to the
  Fourier transform of the basic Euclidean scattering problem,
  $\Delta-\lambda$, $\lambda>0$, which is multiplication by
  $|\xi|^2-\lambda$. The dynamics is exactly as shown above when the
  zero set is a circle.}
\label{fig:fredholm-glob}
\end{figure}

In order to achieve this, in the simplest setting of compact manifolds
without boundary, $X$, we require a
non-trapping dynamics for the formally self-adjoint operator $P$ of
order $m$. Here non-trapping is understood in
the sense that the characteristic set $\Sigma$ of $P$ has connected
components $\Sigma_j$, $j=1,\ldots,k$, in each of which one is given
smooth conic submanifolds $\Lambda_{j,\pm}$ (with
$\Lambda_\pm=\cup_j\Lambda_{j,\pm}$) which act as normal sources
($-$) or sinks ($+$) for the bicharacteristic flow within $\Sigma_j$
in a precise sense described above,
and all bicharacteristics in $\Sigma_j$ except those in
$\Lambda_{j,\pm}$, tend to $\Lambda_{j,+}$ in the forward and to
$\Lambda_{j,-}$ in the backward direction (relative to the flow
parameter) along the bicharacteristic flow, see Figure~\ref{fig:fredholm-glob}. (As recalled above, this setup
can be generalized further, for instance it is stable under general
perturbations in $\Psi^m(X)$ even though the details of the dynamics
are not such in general.) In this case, on variable
order weighted Sobolev spaces $H^s$, with $s$ monotone
increasing/decreasing in each component of the characteristic set
along the Hamilton flow, and satisfying threshold inequalities at
$\Lambda_{j,\pm}$, $P:\cX\to\cY$ is Fredholm, where
\begin{equation}\label{eq:Fredholm}
\cX=\{u\in H^s:\ Pu\in H^{s-m+1}\},\ \cY=H^{s-m+1}.
\end{equation}
Here the Fredholm estimates take the form
\begin{equation}\begin{aligned}\label{eq:Fredholm-est}
\|u\|_{H^s}&\leq C(\|Pu\|_{H^{s-m+1}}+\|u\|_{H^{r}}),\\
\|v\|_{H^{s'}}&\leq C(\|P^*v\|_{H^{s'-m+1}}+\|v\|_{H^{r'}}),
\end{aligned}\end{equation}
for appropriate $r,r'$ with compact inclusion $H^s\to H^r, H^{s'}\to
H^{r'}$, where we take $s'=-s+m-1$, so $s'-m+1=-s$. Note that with
this choice of $s'$ the space on the left hand side, resp.\ in the
first term on the right hand side, of the first inequality is the dual (relative to $L^2$) of
the first space of the right hand side, resp.\ the left hand side of
the second inequality, as required for the functional analytic setup.
Here \eqref{eq:Fredholm-est} is an estimate in terms of Sobolev spaces (which $\cY$ is, but
$\cX$ is not), but it implies the
Fredholm property \eqref{eq:Fredholm}; see \cite[Section~2.6]{Vasy-Dyatlov:Microlocal-Kerr}.

If $P=P^*$, then the threshold regularity is $(m-1)/2$, i.e.\ $s$ can
be almost constant, but it has to be slightly below $(m-1)/2$ at one
end of each bicharacteristic, and slightly above $(m-1)/2$ at the other.
{\em Assuming that these problems are invertible}, the inverse is
independent of the choice of $s$ in a natural sense, as long as the
increasing/decreasing direction of $s$ is kept unchanged along each
component of the characteristic set (see \cite[Remark~2.9]{Vasy-Dyatlov:Microlocal-Kerr}). Note that
in the case of 
invertibility, the compact term of \eqref{eq:Fredholm-est} can be dropped, and one concludes that 
$P^{-1}:H^{s-m+1}\to H^s$, $(P^*)^{-1}:H^{-s}\to H^{-s+m-1}$ are
bounded maps, with 
$(P^*)^{-1}=(P^{-1})^*$. 
(Here invertibility is not a
serious issue for our purposes; see Remark~\ref{rem:non-inv}.) Letting $I\subset\{1,\ldots,k\}=J_k$
be the subset on which $s$ is increasing (i.e.\ where estimates are
propagated {\em backwards}), we denote by
$$
P^{-1}_I:\cY_I\to\cX_I
$$
the
corresponding inverse; here $\cX_I,\cY_I$ stand for the spaces
$\cX,\cY$ above for any choice of $s$ compatible with $I$. Thus, $P_\emptyset^{-1}$ is the Feynman, or forward
propagator, i.e.\ it propagates estimates $H_p$-forward along the
bicharacteristics, so for $\phi\in \CI(X)$,
$\WF(P_\emptyset^{-1}\phi)\subset \cup_j\Lambda_{j,+}$, while
$P_{J_k}^{-1}$ is the backward, or anti-Feynman, propagator. For
general $\phi\in H^{s-m+1}$, $\WF(P^{-1}_\emptyset \phi)$ is contained in the
image of $\WF(\phi)$ under the forward Hamilton flow (interpreted so
that the image of the sources under the forward flow is all
bicharacteristics emanating from them) union the sinks $\cup_j\Lambda_{j,+}$; the analogous
statement for the backward flow holds for $\WF(P^{-1}_{J_k} \phi)$. Such a setup is explained
in detail in \cite[Section~2]{Vasy-Dyatlov:Microlocal-Kerr}.

We recall an example, which will also be used below, from
\cite{Vasy-Dyatlov:Microlocal-Kerr}, given in this form in
\cite{Vasy:Resolvents-AH-dS}. If one considers the Minkowski wave
operator $\Box_g$ on $\RR^{n+1}_{z,t}$, or more conveniently
$x^{-(n-2)/2-2}\Box_g x^{(n-2)/2}$, with $x=(|z|^2+t^2)^{-1/2}$, then
the Mellin transform of this operator in the radial variable on
$\RR^{n+1}$, or its reciprocal $x$, is a family of operators on the
sphere $\sphere^n$; here $\sphere^n$ arises as a smooth transversal to the dilation
orbits on $\RR^{n+1}\setminus o$. This family $P_\sigma$, depending on the Mellin
dual parameter $\sigma$, is an example of this setup with
$X=\sphere^n$. As explained in \cite{Vasy:Resolvents-AH-dS}, in fact
$P_\sigma$ is elliptic/hyperbolic in the region of $\sphere^{n}$ interior/exterior of the
Minkowski light cone; it turns out to be related to the spectral
family of the Laplacian on hyperbolic space, resp.\ the d'Alembertian
on de Sitter space. This example, and natural generalizations, such as
the spectral family of Laplacian and the Klein-Gordon operator on
{\em even asymptotically hyperbolic and de Sitter spaces} (even on
differential forms),
respectively, discussed
in \cite{Vasy:Microlocal-AH} and \cite{Vasy:Analytic-forms}, will arise
again later in this paper.

The basic idea of such a compact dynamical setup first appeared in
Melrose's work on scattering \cite{RBMSpec}, where $P=\Delta-\lambda$,
$\Delta$ is the Laplacian of a scattering metric (large end of a
cone), $\lambda>0$. In that case there are only two propagators, whose
difference is essentially the spectral measure, so the positivity
statement is that of the spectral measure for $\Delta$. In some sense
then, while one should not think of $P$ as a self-adjoint operator in
our general setting (though it is formally self-adjoint), since the
adjoint propagates estimates always the opposite way (corresponding to
having to work in dual function spaces), one still has a positivity
property analogous to these spectral measures. Indeed, from a certain
perspective, the proof given below is inspired by an analogous proof
in scattering theory, related to Melrose's `boundary pairing' \cite[Section~13]{RBMSpec}, though in that setting there are more standard
proofs as well. We refer to the discussion around
Theorem~\ref{thm:sc-elliptic} for more detail.

The main result is the following.

\begin{thm}\label{thm:main}
Suppose $P=P^*\in\Psi(X)$ is as above (i.e.\ $X$ is compact, the
principal symbol $p$ is real, the Hamilton dynamics is non-trapping),
possibly acting on a vector bundle with scalar principal symbol. If
$P^{-1}_{J_k}$, $P^{-1}_\emptyset$ exist (rather than $P$ being merely
Fredholm between the appropriate spaces)
then the operator $\imath(P^{-1}_{J_k}-P^{-1}_\emptyset)$ is
positive, i.e.\ it is symmetric
\begin{equation}\label{eq:symmetry}
\langle \imath(P^{-1}_{J_k}-P^{-1}_\emptyset)\phi,\psi\rangle=\langle
\phi,\imath(P^{-1}_{J_k}-P^{-1}_\emptyset)\psi\rangle,\qquad \phi,\psi\in\CI(X),
\end{equation}
and for all $\phi\in\CI(X)$,
\begin{equation}\label{eq:positive}
\langle \imath(P^{-1}_{J_k}-P^{-1}_\emptyset)\phi,\phi\rangle\geq 0.
\end{equation}
\end{thm}

\begin{rem}
Note that this is the same formula (in particular the sign matches) as in
\cite[Theorem~6.6.2]{FIOII}, with $S_{\tilde n}$ replaced by
$E_{\tilde N}-E_\emptyset$ in \cite[Theorem~6.6.2]{FIOII}, for there $S_{\tilde n}$ is relative to
the backward propagator, denoted by $E_\emptyset$ there.

Note also that the proof given below also shows the symmetry of
$\imath(P^{-1}_{I^c}-P^{-1}_I)$ for any $I\subset J_k$, assuming these
inverses exist (rather than $P$ being just Fredholm between the
corresponding spaces) although
positivity properties are lost. However, see Corollary~\ref{cor:main}
for a parametrix statement, and Remark~\ref{rem:non-inv} regarding invertibility.
\end{rem}

\begin{rem}\label{rem:non-inv}
As the following proof shows, only minor changes are needed if $P$ is
merely Fredholm between the appropriate spaces. Namely for each $I$
let $\cW_I$ be a complementary subspace to the finite dimensional
subspace $\Ker_I P=\Ker P$ of $\cX_I$. Then for $\phi\in\Ran_I P=\Ran
P\subset\cY_I$ there exists a unique $u\in\cW_I$ such that $Pu=\phi$;
we may define $P_I^{-1}\phi=u$. Then \eqref{eq:symmetry} holds if we
require in addition
$\phi,\psi\in\Ran_\emptyset P\cap\Ran_{J_k} P$ and \eqref{eq:positive}
holds if we require $\phi\in\Ran_\emptyset P\cap\Ran_{J_k} P$, as follows
immediately from the proof we give below. (These are finite
codimension conditions!) Note that different choices
of $\cW_I$ do not affect either of the inner products
\eqref{eq:symmetry}-\eqref{eq:positive} since $P_I^{-1}$ applied to an
element $\phi$ of $\Ran_IP$ is being paired with an element of $\Ran_{I^c}P$,
and the latter annihilates (i.e.\ is orthogonal with respect to the
$L^2$ pairing) $\Ker_I P$, i.e.\ changing $P_I^{-1}\phi$ by an element
of $\Ker_I P$ leaves the inner product unchanged.

We now discuss what happens under an additional hypothesis,
$\Ker_\emptyset P,\Ker_{J_k}P\subset\CI(X)$. In this case,
$\Ker_\emptyset P=\Ker_{J_k} P$ (since elements of both are simply
elements of $\CI(X)$ which are mapped to $0$ by $P$); denote this
finite dimensional
space by $\cF$. In this case one
can use the $L^2$-orthocomplements of $\cF$ to define $\cW_\emptyset$
and $\cW_{J_k}$ in $\cX_\emptyset$, resp.\ $\cX_{J_k}$. That is,
$\cW_\emptyset$, resp.\ $\cW_{J_k}$ are the subspaces of
$\cX_\emptyset$, resp.\ $\cX_{J_k}$, consisting of distributions
$L^2$-orthogonal to $\cF$ (which is a subset of both of these spaces!); this
makes sense since $\cF$ is a subspace of
$\CI(X)$. Similarly, $\cF$ gives orthocomplements to $\Ran_\emptyset
P$ and $\Ran_{J_k} P$ in $\cY_\emptyset$, resp,\ $\cY_{J_k}$.
Thus, one can define $P_I^{-1}\phi$, $\phi\in\cY_I$,
$I=\emptyset,J_k$, by defining it to be $P_I^{-1}\phi_1$,
$\phi=\phi_1+\phi_2\in\Ran_IP+\cF$,
where $P_I^{-1}$ takes values in $\cW_I$. In this case, the inner products
\eqref{eq:symmetry}-\eqref{eq:positive} are unaffected by the second
component of the function in both slots, and thus they
remain true for all $\phi,\psi\in\CI(X)$.

While the assumption $\Ker_\emptyset P,\Ker_{J_k}P\subset\CI(X)$ may
seem unnatural, one expects it to hold in analogy with scattering
theory: there incoming or outgoing elements of the tempered distributional
nullspace of the operator necessarily vanish, giving that any such
element is necessarily Schwartz. This can in fact be proved in the
present setting as well using the functional calculus for an elliptic
operator; since this is a bit involved we defer this to another paper,
and we choose to discuss this here only in
the setting of Theorem~\ref{thm:sc-elliptic} for differential operators,
where this is straightforward since the role of the elliptic operator
is played by the weight $x$.
\end{rem}

\begin{proof}
The symmetry statement is standard; one can arrange the function
spaces so that $P^{-1}_{J_k}$ is exactly the inverse of $P^*=P$ on
(essentially) the
duals of the spaces (in reversed role) on which $P_\emptyset^{-1}$
inverts $P$, so $P^{-1}_{J_k}=(P^{-1}_\emptyset)^*$, see
\cite[Section~2]{Vasy-Dyatlov:Microlocal-Kerr} and
\eqref{eq:Fredholm-est} above. Here `essentially'
refers to the fact that the Fredholm estimates \eqref{eq:Fredholm-est},
with the compact terms dropped as remarked above, due to
invertibility, give bounded maps
$P^{-1}:H^{s-m+1}\to H^s$, $(P^*)^{-1}:H^{-s}\to H^{-s+m-1}$, with
$(P^*)^{-1}=(P^{-1})^*$. Correspondingly, the symmetry actually holds
for any $\phi,\psi\in H^{-s}\cap H^{s-m+1}$, $s$ satisfying the requirements for the $\emptyset$-inverse.

We turn to the proof of positivity, with $I=J_k$ to minimize
double subscripts. Let $\smooth_r$, $r\in(0,1)$
be a family of (finitely) smoothing operators, converging to $\smooth_0=\Id$ as
$r\to 0$ in the usual manner, so $\smooth_r\in\Psi^{-N}(X)$, $N>1$ for
$r\in(0,1)$, $\smooth_r$, $r\in(0,1)$ is uniformly bounded in
$\Psi^0(X)$, converging to $\Id$ in $\Psi^\ep(X)$ for all $\ep>0$. Concretely, with $\rho$ a defining
function of $S^*X=\pa\overline{T}^*X$ (e.g.\ homogeneous of degree
$-1$ away from the zero section), we can let the principal symbol
$\symb_r$ of $\smooth_r$ be $(1+r\rho^{-1})^{-N}$, $N>1$. Let $u_I=P_I^{-1}\phi$,
$u_\emptyset=P_\emptyset^{-1}\phi$. Then for $\phi\in\CI(X)$, as $Pu_I=\phi=Pu_\emptyset$,
\begin{equation*}\begin{aligned}
&\langle \imath(P^{-1}_I-P^{-1}_\emptyset)\phi,\phi\rangle=\langle
\imath (u_I-u_\emptyset),Pu_\emptyset\rangle\\
&\qquad=\lim_{r\to 0}\langle \imath
\smooth_r(u_I-u_\emptyset),Pu_\emptyset\rangle
=\lim_{r\to 0}\langle \imath [P,\smooth_r](u_I-u_\emptyset),u_\emptyset\rangle.
\end{aligned}\end{equation*}
Now note that $[P,\smooth_r]$ is uniformly bounded in $\Psi^{m-1}(X)$,
converging to $[P,\Id]=0$ in $\Psi^{m-1+\ep}(X)$, $\ep>0$, so $[P,\smooth_r]\to 0$ strongly as a
bounded operator $H^{\sigma}\to H^{\sigma-m+1}$.
By a
standard microlocal argument about distributions with disjoint wave
front set, using also the above statements on $[P,\smooth_r]$, we have
\begin{equation}\label{eq:cross-term}
\lim_{r\to 0}\langle \imath
[P,\smooth_r]u_I,u_\emptyset\rangle=0.
\end{equation}
To see this claim, let $\Lambda_\pm=\cup_j\Lambda_{\pm,j}$,
$Q_+\in\Psi^0(X)$ be such that $\WF'(\Id-Q_+)\cap \Lambda_+=\emptyset$,
$\WF'(Q_+)\cap\Lambda_-=\emptyset$, i.e.\ $Q_+$ is microlocally the
identity at $\Lambda_+$, microlocally $0$ at $\Lambda_-$, we have (as
$I=J_k$) $Q_+u_I\in \CI(X)$, $(\Id-Q_+)u_\emptyset\in\CI(X)$ since
$\WF(u_\emptyset)\subset\Lambda_+$, $\WF(u_I)\subset\Lambda_-$. Define
$Q_-$ similarly, with $\Lambda_\pm$ interchanged, and such that
$\WF'(\Id-Q_+)\cap\WF'(\Id-Q_-)=\emptyset$ (so at each point at least
one of $Q_\pm$ is microlocally the identity); then $Q_- u_\emptyset\in\CI(X)$
and $(\Id-Q_-) u_I\in\CI(X)$. Thus,
\begin{equation*}\begin{aligned}
\langle \imath
[P,\smooth_r]u_I,u_\emptyset\rangle=&\langle \imath
[P,\smooth_r]Q_+ u_I,u_\emptyset\rangle+\langle
Q_-^*[P,\smooth_r](\Id-Q_+) u_I,u_\emptyset\rangle\\
&\qquad+\langle (\Id-Q_-^*)[P,\smooth_r](\Id-Q_+) u_I,u_\emptyset\rangle.
\end{aligned}\end{equation*}
Now the first term goes to $0$ as $r\to 0$ since $Q_+u_I\in\CI(X)$ in
view of the stated strong convergence of the commutator to $0$, while
the second term goes to $0$ similarly due to $Q_-
u_\emptyset\in\CI(X)$ and the stated strong convergence of the
commutator to $0$. Finally, in view of the disjoint wave front set of
$\Id-Q_+$ and $\Id-Q_-$, thus $\Id-Q_+$ and $\Id-Q_-^*$, $
(\Id-Q_-^*)[P,\smooth_r](\Id-Q_+)$ is in fact uniformly bounded in
$\Psi^{-k}(X)$ for any $k$, and indeed converges to $0$ in
$\Psi^{-k}(X)$, so the third term also goes to $0$. This proves \eqref{eq:cross-term}.

So it remains to consider
$-\langle \imath [P,\smooth_r]u_\emptyset,u_\emptyset\rangle$. But the
principal symbol of $\imath[P,\smooth_r]$ is
$$
H_p\symb_r=Nr\rho^{-2}(1+r\rho^{-1})^{-1}(H_p\rho)\symb_r
$$
which takes the form
$c_r^2\symb_r$ at the sources (as $H_p\rho=\rho^{-m+2}\beta_0$ with $\beta_0$ positive there), and
$-c_r^2\symb_r$ at the sinks. In our case, the wave front set of
$u_\emptyset$ is at the sinks $\Lambda_+$, so are concerned about the
second region. Let $c_r$ be a symbol with square
$-Nr\rho^{-2}(1+r\rho^{-1})^{-1}(H_p\rho)\chi_+^2$, where $\chi_+$ is
a cutoff function, identically $1$ near $\Lambda_+$, supported close
to $\Lambda_+$, and letting $C_r$ be a quantization of this with the
quantization arranged using local coordinates and a partition of
unity; these are being specified so that $C_r$ is uniformly bounded in
$\Psi^{(m-1)/2}(X)$, and {\em still tends to $0$ in
  $\Psi^{(m-1)/2+\ep}(X)$ for $\ep>0$}.
Similarly,
let $E_r$ be a quantization of
$Nr\rho^{-2}(1+r\rho^{-1})^{-1}(H_p\rho)(1-\chi_+^2)$. Then we have
$$
\imath[P,\smooth_r]=-C_r^*\tilde\smooth_r^*\tilde\smooth_r C_r+E_r+F_r,
$$
where the principal symbol of $\tilde\smooth_r$ is the square root of
that of $\smooth_r$, where the family $E_r$ is uniformly bounded in
$\Psi^{m-1}(X)$, has (uniform!) wave front set
disjoint from $\Lambda_+$, while $F_r$ is uniformly bounded in
$\Psi^{m-2}(X)$, and further both $E_r$ and $F_r$ tend to $0$ in higher
order pseudodifferential operators. The disjointness of the uniform wave front set of $E_r$
from $\Lambda_+$, thus from the wave front set of
$u_\emptyset$, and further that it tends to $0$ as $r\to 0$ in the relevant
sense, shows by an argument similar to the proof of
\eqref{eq:cross-term} that
$$
\lim_{r\to 0}\langle E_r u_\emptyset,u_\emptyset\rangle=0.
$$
On the other hand, as $u_\emptyset$ is in $H^{(m-1)/2-\ep}$ for all $\ep>0$,
the fact that $F_r\to 0$ in $\Psi^{m-2+\ep}$, $\ep>0$, and thus
$F_r\to 0$ strongly as a family of operators $H^{(m-1)/2-\ep}\to
H^{-(m-1)/2+\ep}$, $\ep<1/2$,  yields that
$$
\lim_{r\to 0}\langle F_r u_\emptyset,u_\emptyset\rangle=0.
$$
Finally,
$$
\langle C_r^*\tilde\smooth_r^*\tilde\smooth_r C_r u_\emptyset,
u_\emptyset\rangle\geq 0
$$
for all $r$, so
\begin{equation*}\begin{aligned}
\langle \imath(P^{-1}_I-P^{-1}_\emptyset)\phi,\phi\rangle&=
\lim_{r\to 0}\langle \imath
[P,\smooth_r](u_I-u_\emptyset),u_\emptyset\rangle\\
&=\lim_{r\to 0}\langle C_r^*\tilde\smooth_r^*\tilde\smooth_r C_r u_\emptyset,
u_\emptyset\rangle\geq 0.
\end{aligned}\end{equation*}
This
proves the theorem.
\end{proof}

Before proceeding, we now discuss generalized inverses for $P_I$ when
$P_I$ is not invertible, rather merely Fredholm. Note that since
$\CI(X)$ is dense in $\cY_I=H^{s_I-m+1}$, the closed subspace of
finite codimension $\Ran_I P$ has a complementary subspace
$\cZ_I\subset\CI(X)$ in $\cY_I$: indeed, the orthocomplement of
$\Ran_I P$ in the Hilbert space $\cY_I$ is finite dimensional, and
approximating an orthonormal basis for it by elements of $\CI(X)$
gives the desired complementary space. We now can decompose an
arbitrary $f\in\cY_I$ as $f=f_1+f_2$, $f_1\in \Ran_I P$, $f_2\in\cZ_I$
and then, letting $\cW_I$ be a complementary subspace of $\cX_I$ to
$\Ker_I P$,
 $f_1=P u$
 for a unique $u\in \cW_I$; we let $P_I^{-1}f=u$, so $P_I^{-1}$ is a
 {\em generalized inverse} for $P$. Note that as
 $\cZ_I\subset\CI(X)$, $\WF^\sigma(f_1)=\WF^\sigma(f)$ for all
 $\sigma$. The propagation of singularities, for $f\in H^\sigma$, $\sigma>-(m-1)/2$,
 $Pu=f_1$, $u\in \cX_I$ shows that $\WF^{\sigma+m-1}(u)\subset\cup_{j\in
   I}\Lambda_{j,-}\cup\cup_{j\in I^c}\Lambda_{j,+}$. This suffices for
 all the arguments below.

An immediate corollary of Theorem~\ref{thm:main} is the Duistermaat-H\"ormander theorem:

\begin{cor}\label{cor:main}(cf.\ Duistermaat and H\"ormander
  \cite[Theorem~6.6.2]{FIOII})
Suppose that $P$ is as in Theorem~\ref{thm:main} (in particular,
$P_\emptyset$, $P_{J_k}$ are invertible).
For all $I$, there exists an operator $\tilde S_I$ such that
$P_I^{-1}-P_\emptyset^{-1}$ differs from $\tilde S_I$ by an operator
that is smoothing away from $\Lambda_\pm$ in the sense that $\phi\in\dist(X)$,
$\WF^\sigma(\phi)\cap (\Lambda_+\cup\Lambda_-)=\emptyset$, $\sigma>-(m-1)/2$, implies
$\WF((P_I^{-1}-P_\emptyset^{-1}-\tilde
S_I)\phi)\subset\Lambda_+\cup\Lambda_-$,
and such that $\tilde S_I$
is skew-adjoint and $\imath\tilde S_I$ is positive.

Here, if $P_I$ is
not invertible (i.e.\ is only Fredholm), the statement holds if in
addition $\phi\in\Ran P_I$ in the sense
of Remark~\ref{rem:non-inv}, and more generally for all $\phi$ as above
if $P_I^{-1}$ is a generalized inverse of $P_I$ defined on $\cY_I$
using a complement $\cZ_I$ to $\Ran P_I$ which is a subspace of $\CI$,
as defined above.
\end{cor}

Thus, here smoothing is understood e.g.\ as a statement that for $\phi\in H^{\sigma}(X)$, where
$\sigma>-(m-1)/2$, the operator in question maps to $\CI(X)$,
microlocally away from $\Lambda_\pm$. In fact, as all the operators in
question can naturally be applied to distributions with wave
front set
away from $\Lambda_\pm$ (by suitable choice of the order function $s$), which is the context of the
Duistermaat-H\"ormander result, and smoothing holds in this extended
context as well, as stated in the corollary.

\begin{rem}
In the case of $P_\emptyset$ and $P_{J_k}$, in 
Remark~\ref{rem:non-inv} we showed that if $\Ker_\emptyset 
P,\Ker_{J_k} P\subset\CI(X)$, then we have canonical generalized 
inverses $P_\emptyset^{-1}$, $P_{J_k}^{-1}$ which satisfy the
properties \eqref{eq:symmetry}-\eqref{eq:positive}. Thus, relaxing the
invertibility hypothesis for $P_\emptyset$, $P_{J_k}$, but under this
additional assumption on the kernels of these operators, conclusion of
this Corollary still holds.
\end{rem}

\begin{proof}
In the following discussion we assume that $P_I$ is invertible. In
fact, all
the arguments go through for a generalized inverse as in the statement
of the theorem, but it is more convenient to not have to write out repeatedly
decompositions with respect to which the generalized inverse is taken.

We use a microlocal partition of unity
$\sum_{j=0}^kB_j$, $B_j=B_j^*$,
with $B_0$ having wave front set in the elliptic set, $B_j$, $j\geq 1$
having wave front set disjoint from the components $\Sigma_l$, $l\neq
j$, of the characteristic set. Let
$$
T_j=B_j
(P_{J_k}^{-1}-P_\emptyset^{-1})B_j.
$$
Then for any $I$,
$$
\tilde S_I=\sum_{j\in I} T_j
$$
has the required properties, with skew-adjointness of $\tilde S_I$ and
positivity of $\imath\tilde S_I$ following from the main theorem
above.

To see the parametrix property, note that for $j\neq 0$, $B_j=\Id$
microlocally near $\Sigma_j$, while $B_l=0$ microlocally near
$\Sigma_j$ for $l\neq j$. Thus, for $\phi\in H^{\sigma}(X)$, where
$\sigma>-(m-1)/2$,
$$
P(P_\emptyset^{-1}+\tilde
S_I)\phi=\phi+\sum_{j\in I}[P,B_j](P_{J_k}^{-1}-P_\emptyset^{-1})B_j\phi,
$$
with the wave front set of the commutator, and thus of all but the
first term, being in the elliptic set. But
$P(P_{J_k}^{-1}-P_\emptyset^{-1})B_j\phi=0$, so microlocal elliptic
regularity shows that $[P,B_j]
(P_{J_k}^{-1}-P_\emptyset^{-1})B_j\phi\in\CI(X)$.

Notice that microlocal elliptic regularity also shows that all
parametrices are microlocally the same in the elliptic set: if
$Pu-Pv\in\CI(X)$, then $u-v$ has wave front set disjoint from the
elliptic set of $P$. So in order to analyze our parametrix, it
suffices to consider the characteristic set.

Microlocally near $\Sigma_j$, $P_\emptyset^{-1}f$,
resp.\ $P_{J_k}^{-1}f$, $f\in  
H^{\sigma}(X)$, 
solve $Pu-f\in\CI(X)$, with
$\WF^{\sigma+m-1}(u)\subset \Lambda_{+,j}$, resp.\
$\WF^{\sigma+m-1}(u)\subset \Lambda_{-,j}$.
Further $P_I^{-1}f$ has the same property as one of these, depending on whether
$j\notin I$ or $j\in I$. In particular, for $j\notin I$,
$u=P_I^{-1}\phi-P_\emptyset^{-1}\phi$ solves $Pu=0$, with
$\WF^{\sigma+m-1}(u)\cap\Sigma_j\subset \Lambda_{+,j}$, which implies
by propagation of singularities (including the version at the radial
points in $\Lambda_{-,j}$, where $u$ is a priori in a better space than
the threshold Sobolev regularity) that in fact
$\WF(u)\cap\Sigma_j\subset \Lambda_{+,j}$.
Since  microlocally near $\Sigma_j$,
$(P_\emptyset^{-1}+\tilde S_I)\phi$ is the same as
$P_\emptyset^{-1}\phi$ if $j\notin I$, we deduce that
$P_I^{-1}-(P_\emptyset^{-1}+\tilde S_I)$ is smoothing near such $j$,
in the sense that in this neighborhood of $\Sigma_j$,
$\WF(P_I^{-1}\phi-(P_\emptyset^{-1}+\tilde S_I)\phi)$ is contained in $\Lambda_{+,j}$, so we only
need to consider $j\in I$.

Since $B_j\phi$ and $\phi$ are the same
microlocally near $\Sigma_j$, by the propagation of singularities,
again using the a priori better than threshold Sobolev regularity at
$\Lambda_{+,j}$, $u=P_{J_k}^{-1}(\phi-B_j\phi)$ has
$\WF(u)\cap\Sigma_j\subset \Lambda_{-,j}$, and similarly for
$P_{\emptyset}^{-1}(\phi-B_j\phi)$ (for $\Lambda_{+,j}$). In view of $B_j$ being
microlocally the identity near $\Sigma_j$, and trivial near
$\Sigma_k$, $k\neq j$, we deduce that the intersection of the wave
front set of $P_{\emptyset}^{-1}
\phi-B_jP_\emptyset^{-1}B_j \phi$ with $\Sigma_j$ is in
$\Lambda_{+,j}$. Similar arguments give that for $j\in I$ the intersection of the wave
front set of $P_{I}^{-1}
\phi-B_jP_{J_k}^{-1}B_j \phi$ with $\Sigma_j$ is in
$\Lambda_{-,j}$. The conclusion is that, microlocally near $\Sigma_j$,
$j\in I$, the wave front set of $(P_\emptyset^{-1}+\tilde
S_I)\phi-P_I^{-1}\phi$ is in $\Lambda_{+,j}\cup\Lambda_{-,j}$.
This proves that
$P_I^{-1}-P_\emptyset^{-1}$ differs from $\tilde S_I$ by an operator
that is smoothing away from $\Lambda_\pm$, completing the proof of the
corollary.
\end{proof}

Notice that while it is a distinguished parametrix in the Duistermaat-H\"ormander
sense, $P_\emptyset^{-1}+\tilde S_I$ is in principle not necessarily one of
our distinguished inverses, $P_I^{-1}$. Indeed, while $P_I^{-1}$ maps
$\phi\in\CI(X)$ to have wave front set disjoint from $\Lambda_{j,+}$
for $j\in I$, on the other hand, for $j\in I$ the difference of
$P_\emptyset^{-1} \phi$ and $B_j P_\emptyset^{-1}B_j\phi$ at
$\Lambda_{j,+}$ is not necessarily smooth, though it does have wave
front set (locally) contained in $\Lambda_{j,+}$ (i.e.\ the difference
is smoothing away from $\Lambda_{j,+}$ within $\Sigma_j$). If $B_j$
can be arranged to commute with $P$, however, this statement can be
improved.

\section{Positivity in Melrose's b-pseudodifferential algebra}
There are natural extensions to b- and scattering settings of Melrose
(see \cite{Melrose:Atiyah} for a general treatment of the b-setting,
\cite{RBMSpec} for the scattering setting), such as
the wave equation and the Klein-Gordon equation on asymptotically
Minkowski spaces, in the sense of `Lorentzian scattering metrics' of
Baskin, Vasy and Wunsch, see
\cite{Baskin-Vasy-Wunsch:Radiation} and \cite[Section~5]{Hintz-Vasy:Semilinear}. This in particular includes the
physically relevant example of Minkowski space (and perturbations of
an appropriate type) that motivated this part of the
Duistermaat-H\"ormander work. Since no new analytic work is necessary
in these new settings (i.e.\ one essentially verbatim repeats the
proof of Theorem~\ref{thm:main} and Corollary~\ref{cor:main}, changing various bundles, etc.), we
only briefly recall the setups and state the corresponding theorems,
explaining any (minor) new features.

Before proceeding, we recall that Melrose's b-analysis is induced by the analysis
of totally characteristic, or b-, differential operators, i.e.\ ones
generated (over $\CI(M)$, as finite sums of products) by vector fields
$V\in\Vb(M)$ tangent to the boundary of a manifold with boundary
$M$. Locally near some point in $X=\pa M$, with the boundary defined by
a function $x$ (so it vanishes non-degenerately and exactly at $\pa
M$), and with $y_j$,
$j=1,\ldots,n-1$, local coordinates on $X$, extended to
$M$, these vector fields are linear combinations of the vector fields
$x\pa_x$ and $\pa_{y_j}$ with smooth
coefficients, i.e.\ are of the form $a_0(x\pa_x)+\sum a_j\pa_{y_j}$. Correspondingly, they are exactly the set of all smooth sections
of a vector bundle, $\Tb M$. Thus, the dual bundle $\Tb^*M$ has smooth
sections locally of the form $b_0\,\frac{dx}{x}+\sum b_j\,dy_j$, with
$b_j$ smooth. Then (classical) b-pseudodifferential operators $P\in\Psib^m(M)$
have principal symbols $p$ which are homogeneous degree $m$ functions
on $\Tb^*M\setminus o$.

Thus, in the b-setting, where this setup was described by
Gell-Redman, Haber and Vasy \cite{Gell-Redman-Haber-Vasy:Feynman},
we require for the strengthened Fredholm framework that $P\in\Psib^m(M)$ is
formally self-adjoint, and the bicharacteristic dynamics in $\Sb^*M$
is as before, i.e.\ with sources and sinks at
$L=L_+\cup L_-\subset\Sb^*M=(\Tb^*M\setminus o)/\RR^+$ (with
$L_+=\Lambda_+/\RR^+$ in the previous notation, where $\Lambda$ was conic). Examples include a modified conjugate of the
Minkowski wave operator, and more generally non-trapping Lorentzian scattering
metrics, namely
if $x$ is a boundary defining function, then the relevant operator is
$P=x^{-(n-2)/2-2}\Box_g x^{(n-2)/2}$ (symmetric with respect to the
b-density $x^{n}\,|dg|$); see \cite{Hintz-Vasy:Semilinear}. The
characteristic set $\Sigma$ satisfies $\Sigma\subset\Sb^*M$, and is a
union of connected components $\Sigma_j$, $j=1,\ldots,k$, just as in the boundaryless
setting. Again, choosing a subset $I$ of $J_k$, we require the
order $s$ to be increasing along the $H_p$-flow on $\Sigma_j$, $j\in I$, decreasing
otherwise, so in $\Sigma_j$, $j\in I$ estimates are propagated
backwards, for $j\in I^c$ forwards. The additional ingredient is to
have a weight $\decay\in\RR$; we then work with the variable order
b-Sobolev spaces $\Hb^{s,\decay}$. The actual numerology of the function
spaces arises from the sources and sinks, namely with $x$ being a
boundary defining function as before and $\rho_\infty$ being a
defining function of fiber infinity in $\overline{\Tb^*}M$ (so e.g.\
can be taken as a homogenous degree $-1$ function on $\Tb^*M$ away from
the zero section), both $H_p x$ and $H_p \rho_\infty$ play a role. A
general numerology is discussed in
\cite[Proposition~2.1]{Hintz-Vasy:Semilinear} for saddle points, with
an analogous numerology also available for other sources/sinks
but is discussed only for $P=x^{-(n-2)/2-2}\Box_g x^{(n-2)/2}$ in \cite[Section~5]{Hintz-Vasy:Semilinear}.
Thus, here for simplicity, we only consider the numerology of
$P=x^{-(n-2)/2-2}\Box_g x^{(n-2)/2}$, though we remark that {\em for
  ultrahyperbolic equations corresponding to quadratic forms on
  $\RR^n$ the numerology is identical}.
The requirement at $L$ then for obtaining the
estimates needed to establish Fredholm properties is $s+\decay>(m-1)/2$
(with $m=2$ for the wave operator) at the components $L_{\pm,j}$ from which
one wants to propagate estimates, and $s+\decay<(m-1)/2$ to which one wants
to propagate estimates. This (plus the required monotonicity of $s$
along bicharacteristics) is still not sufficient, it only gives
estimates of the form
$$
\|u\|_{\Hb^{s,\decay}}\leq C(\|Pu\|_{\Hb^{s-m+1,\decay}}+\|u\|_{\Hb^{\tilde s,\decay}}),
$$
with $\tilde s<s$; here the problem is that the inclusion
$\Hb^{s,\decay}\to\Hb^{\tilde s,\decay}$ is {\em not} compact, because there is
no gain in the decay order, $\decay$.
Thus, one needs an
additional condition involving the Mellin transformed normal operator,
$\hat P(.)$.

One arrives at the normal operator by `freezing
coefficients' at $X=\pa M$, namely by using a collar neighborhood
$X\times[0,\ep)_x$ of
$X$, including it in $X\times[0,\infty)$, obtaining an operator by
evaluating the coefficients of $P$ at $x=0$ (which can be done in a
natural sense) and then regarding the resulting $N(P)$ as a dilation invariant operator
on $X\times\RR^+$, with dilations acting in the second factor. The
Mellin transform then is simply the Mellin transform in the
$\RR^+$-factor. Thus, the Mellin transformed normal operator is a family of operators, $\Cx\ni\sigma\mapsto\hat
P(\sigma)$, on $X=\pa M$. In fact, this is an analytic Fredholm family by the
boundaryless analysis explained above (with the dynamical assumptions
on $P$ implying those for $\hat P(\sigma)$), which in addition has the
property that for any $C>0$ it is invertible
in $|\im\sigma|<C$ for $|\sigma|$ large (with `large' depending on
$C$), due to the high energy, or semiclassical version, of these
Fredholm estimates. The poles of the inverse are called resonances
and form a discrete set of $\Cx$, with only finitely many in any strip
$|\im\sigma|<C$. If $\decay$ is chosen so that there are no resonances with
$\im\sigma=-\decay$, and if the requirement on $s$ is strengthened to
$s+\decay-1>(m-1)/2$ at the components from which we propagate estimates then $P:\cX\to\cY$ is Fredholm, where
$$
\cX=\{u\in \Hb^{s,\decay}:\ Pu\in \Hb^{s-m+1,\decay}\},\ \cY=\Hb^{s-m+1,\decay}.
$$
(Here the stronger requirement $s+\decay-1>(m-1)/2$ enters when combining
the normal operator estimates with the symbolic estimates, see
\cite[Proposition~2.3 and Section~5]{Hintz-Vasy:Semilinear} and \cite[Theorem~3.3]{Gell-Redman-Haber-Vasy:Feynman}.)
Again, for given $\decay$, if $P$ is actually invertible, $P^{-1}$ only depends on the
choice of $I$ (modulo the natural identification), so we write
$P_I^{-1}$; if we allow $\decay$ to vary it is still independent of $\decay$ as
long as we do not cross any resonances, i.e.\ if $\decay$ and $\decay'$ are such
that there are no resonances $\sigma$ with $-\im\sigma\in[\decay,\decay']$ (if $\decay<\decay'$). Then the arguments given above, with regularization
$\smooth_r$ needed only in the differentiability (not decay) sense, so
$\smooth_r\in\Psib^{-N}(M)$ for $r>0$, uniformly bounded in
$\Psib^0(M)$, converging to $\Id$ in $\Psib^\ep(M)$ for any $\ep>0$
apply if we take the decay order to be $\decay=0$, i.e.\ work with spaces
$\Hb^{s,0}$, the point being that $[P,\smooth_r]\to 0$ in
$\Psib^{m-1+\ep}$ then (there is no extra decay at $X$), so we need to
make sure that $u_I$ lie in a weighted space with weight $0$ to get
the required boundedness and convergence properties. In summary,
this show immediately the following theorem
and corollary, with the analogue of Remark~\ref{rem:non-inv} also valid:

\begin{thm}\label{thm:main-b}
Suppose $P=P^*\in\Psib^m(M)$ is as above, and suppose that no
resonances lie on the real line, $\im\sigma=0$.  If
$P^{-1}_{J_k}$, $P^{-1}_\emptyset$ exist (rather than $P$ being merely
Fredholm between the appropriate spaces)
then the operator $\imath(P^{-1}_{J_k}-P^{-1}_\emptyset)$ is
positive, i.e.\ it is symmetric and for all $\phi\in\dCI(M)$,
$$
\langle \imath(P^{-1}_{J_k}-P^{-1}_\emptyset)\phi,\phi\rangle\geq 0.
$$
\end{thm}

\begin{cor}(cf.\ Duistermaat and H\"ormander
  \cite[Theorem~6.6.2]{FIOII})
Suppose that $P$ is as in Theorem~\ref{thm:main-b} (in particular,
$P_\emptyset$, $P_{J_k}$ are invertible).
For all $I$, there exists an operator $\tilde S_I$ such that
$P_I^{-1}-P_\emptyset^{-1}$ differs from $\tilde S_I$ by an operator
that is smoothing away from $L_\pm$ in the sense that
$\WFb^{\sigma,0}(\phi)\cap (L_+\cup L_-)=\emptyset$,
$\sigma>1-\frac{m-1}{2}$, implies that
$\WFb^{\infty,0}((P_I^{-1}-P_\emptyset^{-1}-\tilde
S_I)\phi)\subset L_+\cup L_-$, and such that $\tilde S_I$
is skew-adjoint and $\imath\tilde S_I$ is positive.

If $P_I$ is not invertible, $P_I^{-1}$ is understood as a generalized
inverse, using a $\dCI(M)$-complement to $\Ran_I P$, similarly to the
discussion preceding Corollary~\ref{cor:main}.
\end{cor}

\section{Positivity in Melrose's scattering pseudodifferential algebra}
The scattering setting, $P\in\Psisc^m(M)$ (one can also have a weight
$l$; this is irrelevant here), is analogous to the b-setting, except that all the principal
symbols are functions (there is no normal operator family), but they
are objects on two intersecting boundary hypersurfaces of the
cotangent bundles: fiber infinity $\Ssc^*M$, and base
infinity $\overline{\Tsc^*}_{\pa M}M$, and (full) ellipticity is the
invertibility of both of these. (Note that these two parts of the
principal symbol agree at the corner
$\Ssc^*_{\pa M}M=\pa\overline{\Tsc^*}_{\pa M}M$ of
$\overline{\Tsc^*}M$.) While here we used the invariant formulation,
an example to which it can {\em always} be locally reduced is the
radial compactification $M=\overline{\RR^n}$ of $\RR^n$; in that case
$\Tsc^*M=\overline{\RR^n}\times\RR^n$ with basis of sections of
$\Tsc^*M$ given by the lifts of the standard coordinate differentials
$dz_j$, $j=1,\ldots,n$, and $\overline{\Tsc^*}M=\overline{\RR^n}\times \overline{\RR^n}$.
This is the setting which Melrose introduced for studying the scattering
theory of asymptotically Euclidean spaces \cite{RBMSpec}; these are compactified
Riemannian manifolds $M$ (so one has a Riemannian metric on $M^\circ$) which are asymptotically the
large ends of cones. For Melrose's problem, the operator $P=\Delta-\lambda$ is elliptic
at fiber infinity, $\Ssc^*M$; note that $\lambda$ is {\em not} lower
order than $\Delta$ in the sense of the relevant principal symbol,
namely at base infinity.

For such scattering problems the previous discussion can
be repeated almost verbatim. Here one works with variable order scattering Sobolev
spaces $\Hsc^{s,\decay}(M)$, with $\ell$ being necessarily variable now due to the
ellipticity at $\Ssc^*M$, see \cite{Vasy:Propagation-Notes}. Again, the relevant dynamical assumption is source/sink
bundles $L_\pm$, where now since we have ellipticity at $\Ssc^*M$, we have $L_\pm\subset
\overline{\Tsc^*}_{\pa M}M$, where now the requirement is $\decay>-1/2$ at the
components from which we want to propagate estimates, and $\decay<-1/2$
at the components
towards which we want to propagate estimates. (For a general operator
of order $l$, the threshold would be $(l-1)/2$, i.e.\ $l$ simply plays
the analogue of the differential order $m$ discussed in the compact
setting $X$.) Actually as above, one
can weaken the assumptions on the dynamics significantly, so one does
not even need a source/sink manifold: one needs a source/sink region,
with suitable behavior in the normal variables. (So for instance, the
more typical lower dimensional sources/sinks/saddles of
\cite{Hassell-Melrose-Vasy:Microlocal} are fine as well for this
analysis; one regards the whole region on the `outgoing' side a sink,
on the `incoming' side a source, regardless of the detailed dynamical behavior.) One then has that
$P:\cX\to\cY$ is Fredholm, where
$$
\cX=\{u\in \Hsc^{s,\decay}:\ Pu\in \Hsc^{s-m,\decay-1}\},\ \cY=\Hsc^{s-m,\decay-1}.
$$
Now our `smoothing' $\smooth_r$ is actually just decay gaining, i.e.\ spatial
regularization, corresponding to $\decay$; this does not affect the proof
of the analogue of the main theorem.
We thus have, with the above notation, with the analogue of
Remark~\ref{rem:non-inv} also holding:

\begin{thm}\label{thm:sc-elliptic}
Suppose $P=P^*\in\Psisc^m(M)$ is as above, in particular elliptic at $\Ssc^*M$.  If
$P^{-1}_{J_k}$, $P^{-1}_\emptyset$ exist (rather than $P$ being merely
Fredholm between the appropriate spaces)
then the operator $\imath(P^{-1}_{J_k}-P^{-1}_\emptyset)$ is
positive, i.e.\ it is symmetric and for all $\phi\in\dCI(M)$,
$$
\langle \imath(P^{-1}_{J_k}-P^{-1}_\emptyset)\phi,\phi\rangle\geq 0.
$$
\end{thm}

\begin{cor}(cf.\ Duistermaat and H\"ormander \cite[Theorem~6.6.2]{FIOII})
Suppose that $P$ is as in Theorem~\ref{thm:sc-elliptic} (in particular,
$P_\emptyset$, $P_{J_k}$ are invertible).
For all $I$, there exists an operator $\tilde S_I$ such that
$P_I^{-1}-P_\emptyset^{-1}$ differs from $\tilde S_I$ by an operator
that is smoothing away from $L_\pm$ in the sense that
$\WFsc^{s-m,\mu}(\phi)\cap (L_+\cup L_-)=\emptyset$, $\mu>1/2$, implies
$\WFsc((P_I^{-1}-P_\emptyset^{-1}-\tilde
S_I)\phi)\subset L_+\cup L_-$, and such that $\tilde S_I$
is skew-adjoint and $\imath\tilde S_I$ is positive.

If $P_I$ is not invertible, $P_I^{-1}$ is understood as a generalized
inverse, using a $\dCI(M)$-complement to $\Ran_I P$, similarly to the
discussion preceding Corollary~\ref{cor:main}.
\end{cor}

Notice that in this setting in fact $P$ is actually self-adjoint on
$L^2_\scl(M)=\Hsc^{0,0}(M)$ as an unbounded operator, which in turn
follows from the invertibility of
$$
P\pm \imath:\Hsc^{s,\decay}\to\Hsc^{s-m,\decay}
$$
for any $s,\decay$; note that $P\pm \imath$ is fully elliptic so invertibility
as a map between any such pair of Sobolev spaces is equivalent to
invertibility between any other pair. In the case of
$P=\Delta_g+V-\lambda$, $g$ as scattering metric, $V\in x\CI(M)$ real, this problem was
studied by Melrose \cite{RBMSpec}, but of course there is extensive
literature in Euclidean scattering theory from much earlier. Then for $\lambda>0$ the
limits
$$
(P\pm \imath 0)^{-1}=\lim_{\ep\to 0}(P\pm \imath\ep)^{-1}
$$
exist in appropriate function
spaces (this is the limiting absorption principle), and
$$
\imath(P^{-1}_{J_k}-P^{-1}_\emptyset)=\imath ((P+\imath
0)^{-1}-(P-\imath 0)^{-1})
$$
is, up to a positive factor, the density
of the spectral measure by Stone's theorem. A direct scattering theory
formula for it, implying its positivity, was given in
\cite[Lemma~5.2]{Hassell-Vasy:Spectral} using the Poisson operators;
this formula in turn arose from `boundary pairings'. This explains in detail
the earlier statement that our result is a generalization of the
positivity of the spectral measure in a natural sense.

This also gives rise to another interesting example, namely an
asymptotically Euclidean space whose boundary has two connected
components, e.g.\ two copies of $\RR^n$ glued in a compact
region. Then the previous theory applies in particular, with the
Feynman and anti-Feynman propagators giving the limiting absorption
principle resolvents. However, one can also work with different
function spaces, propagating estimates forward in one component of the
boundary (and hence the characteristic set), and backward in the
other, relative to the Hamilton flow. The resulting problem is
Fredholm, though the invertibility properties are unclear. This
problem is an analogue of the retarded and advanced propagators (and
thus the Cauchy problem) for the wave equation.

We now discuss the comments in the final paragraph of Remark~\ref{rem:non-inv} in more
detail. For operators of the kind $P=\Delta_g+V-\lambda$, $g$ as
scattering metric, $V\in x\CI(M)$ real, using the
boundary pairing formula Melrose showed in
\cite{RBMSpec} that the nullspaces of $P_\emptyset$ and
$P_{J_k}$ are necessarily in $\dCI(M)$; he then used H\"ormander's
unique continuation theorem to show that in fact these nullspaces are
trivial. There is a more robust proof of these results by a different
commutator approach which, as far as the author knows, goes back to
Isozaki's work in $N$-body scattering \cite[Lemma~4.5]{IsoRad}. In a
geometric $N$-body setting this proof was adapted by Vasy in
\cite[Proposition~17.8]{Vasy:Propagation-2}; it in particular applies
to operators like $P=\Delta_g+V-\lambda$. The argument relies on a
family of commutants given by functions which are {\em not} (uniformly) bounded in the
relevant space of (scattering) pseudodifferential operators, but for
which the commutators themselves {\em are} bounded, and have a
sign modulo lower order terms. In the general setting of pseudodifferential operators, an analogous argument works provided one
uses the functional calculus for an elliptic operator (the weight in
the commutant). One has to be rather careful here because the
commutant family is not bounded: this is the reason that the argument
only implies that elements of the nullspace are in $\dCI(M)$, not that
those with $Pu\in\dCI(M)$ are such; for pairings involving the
commutant and $Pu$ must vanish identically. This point will be
addressed in a future paper in full detail.
Notice that this result {\em only} applies to $P_\emptyset$ and
$P_{J_k}$ as illustrated by the two Euclidean end problem in a
particularly simple setting: the line $\RR_z$ with $V=0$, $\lambda=1$. Then
the complex exponentials $e^{\pm \imath z\cdot\zeta}$ are incoming at
one end, outgoing at the other, thus are in the nullspace of $P_I$,
resp.\ $P_{I^c}$, for $I$ corresponding to the appropriate non-Feynman choice.

The simplest non-elliptic (in
the usual sense) interesting example in the scattering setting is the Klein-Gordon equation on
asymptotically Minkowski like spaces (in the same sense as above, in
the b-case, i.e.\ Lorentzian scattering spaces of \cite{Baskin-Vasy-Wunsch:Radiation}). Here one works with variable order scattering Sobolev
spaces $\Hsc^{s,\decay}(M)$, see \cite{Vasy:Propagation-Notes}. Let
$\rho_\infty$ be a defining function for fiber infinity, $\Ssc^*M$,
and $\rho_{\pa M}$ a defining function for base infinity
$\overline{\Tsc^*}_{\pa M} M$. Again, the relevant dynamical assumption is source/sink
bundles $L_\pm$, where now for simplicity we assume that $L_\pm\subset
\overline{\Tsc^*}_{\pa M}M$ transversal to the boundary of the fiber
compactification and now $\beta_0=\mp \rho_\infty^{m-1}\rho_{\pa M}^{-1}H_p\rho_{\pa M}$ is positive at $L_\pm$
while $\rho_{\infty}^{m-2}H_p\rho_{\infty}$ vanishes
there. In this case, as shown in \cite[Proposition~0.11]{Vasy:Propagation-Notes} (where
the roles of $\rho_\infty$ and $\rho_{\pa M}$ are reversed), the requirement
for propagation estimates at the sources/sinks is $\decay>-1/2$ at the
components from which we want to propagate estimates, and $\decay<-1/2$
at the components
towards which we want to propagate estimates. Actually as above, one
can weaken the assumptions on the dynamics significantly, so one does
not even need a source/sink manifold: one needs a source/sink region,
with suitable behavior in the normal variables. (So for instance, the
more typical lower dimensional sources/sinks/saddles of
\cite{Hassell-Melrose-Vasy:Microlocal} are fine as well for this
analysis; one regards the whole region on the `outgoing' side a sink,
on the `incoming' side a source, regardless of the detailed dynamical
behavior.) With $\decay$ chosen monotone along the $H_p$-flow, satisfying
these inequalities, and with the dynamics being non-trapping in the
same sense as before, one then has that
$P:\cX\to\cY$ is Fredholm, where
$$
\cX=\{u\in \Hsc^{s,\decay}:\ Pu\in \Hsc^{s-m+1,\decay-1}\},\ \cY=\Hsc^{s-m+1,\decay-1}.
$$
Since there are no restrictions on $s$, we may simply take it high
enough so that there are no issues with pairings, etc., as far as $s$
is concerned, and so we do not need to regularize in $s$.
Thus, with the above notation and with the same proof, with
$\smooth_r$ regularizing only in decay:

\begin{thm}\label{thm:main-sc}
Suppose $P=P^*\in\Psisc^m(M)$ is as above.  If
$P^{-1}_{J_k}$, $P^{-1}_\emptyset$ exist (rather than $P$ being merely
Fredholm between the appropriate spaces)
then the operator $\imath(P^{-1}_{J_k}-P^{-1}_\emptyset)$ is
positive, i.e.\ it is symmetric and for all $\phi\in\dCI(M)$,
$$
\langle \imath(P^{-1}_{J_k}-P^{-1}_\emptyset)\phi,\phi\rangle\geq 0.
$$
\end{thm}

\begin{cor}(cf.\ Duistermaat and H\"ormander \cite[Theorem~6.6.2]{FIOII})
Suppose that $P$ is as in Theorem~\ref{thm:main-sc} (in particular,
$P_\emptyset$, $P_{J_k}$ are invertible).
For all $I$, there exists an operator $\tilde S_I$ such that
$P_I^{-1}-P_\emptyset^{-1}$ differs from $\tilde S_I$ by an operator
that is smoothing away from $L_\pm$ in the sense that
$\WFsc^{s-m+1,\mu}(\phi)\cap (L_+\cup L_-)=\emptyset$, $\mu>1/2$, implies
$\WFsc((P_I^{-1}-P_\emptyset^{-1}-\tilde
S_I)\phi)\subset L_+\cup L_-$, and such that $\tilde S_I$
is skew-adjoint and $\imath\tilde S_I$ is positive.

If $P_I$ is not invertible, $P_I^{-1}$ is understood as a generalized
inverse, using a $\dCI(M)$-complement to $\Ran_I P$, similarly to the
discussion preceding Corollary~\ref{cor:main}.
\end{cor}

\section{Asymptotically de Sitter problems}

We end this paper by discussing a new direction.
An interesting class of Lorentzian spaces whose behavior is more
complicated is asymptotically
de Sitter spaces. As shown in \cite{Vasy-Dyatlov:Microlocal-Kerr},
\cite{Vasy:Microlocal-AH} and \cite{Vasy:Resolvents-AH-dS}, the
Klein-Gordon operator $\Box_{X_0}-(n-1)^2/4-\sigma^2$ on these spaces $X_0$
can be analyzed by `capping them off' with asymptotically hyperbolic
spaces $X_\pm$ to obtain a compact manifold without boundary $X$.
(In general, for topological reasons, one needs two copies of the asymptotically de Sitter spaces, see \cite[Section~3]{Vasy:Resolvents-AH-dS}.)
Then on $X$ one has exactly the setup analyzed at the beginning of
this paper. In particular, with the characteristic set having two
components (if only a single connected asymptotically de Sitter space
was used) one has forward and backward propagators, which propagate
estimates in the opposite direction relative to the Hamilton flow
in the two components, as
well as Feynman and anti-Feynman propagators which propagate either
forward everywhere along the Hamilton flow or backward everywhere. In
the aforementioned papers the connection between the forward and
backward propagators on $X$ and the resolvents of the Laplacian on
$X_\pm$ as well as the forward and backward propagators on $X_0$ is
explained; see in particular
\cite[Section~4]{Vasy:Resolvents-AH-dS}. For instance, if $F=\{j\}$,
where $\Sigma_j$ is the component of $\Sigma$ on which the de Sitter
time function is decreasing along the bicharacteristics, then
$P_F^{-1}$ gives rise to the forward propagator
$$
(\Box_{X_0}-(n-1)^2/4-\sigma^2)^{-1}_\future=x_{X_0}^{-\imath\sigma+(n-1)/2}P_F^{-1}x_{X_0}^{\imath\sigma-(n-1)/2-2},
$$
where $x_{X_0}$ is a boundary defining function of $X_0$ (which is
thus time-like near $\pa X_0$).
In particular, these global
propagators on $X$ can be used to analyze the local objects on $X_0$
and $X_\pm$; this is essentially a consequence of the evolution
equation nature of the wave equation in the de Sitter region. Thus,
for instance, it does not matter how one caps off $X_0$ above, the
forward propagator on $X$, in the appropriate sense (conjugation and
multiplication) restricts to the forward propagator on $X_0$ --- an
object independent of the choice of the caps $X_\pm$!

A natural question is then whether this method allows one to define a
canonical Feynman propagator on $X_0$. Certainly one choice arises by
taking $P_\emptyset^{-1}$ on $X$, and letting
$$
(\Box_{X_0}-(n-1)^2/4-\sigma^2)^{-1}_\Feynman=x_{X_0}^{-\imath\sigma+(n-1)/2}P_\emptyset^{-1}x_{X_0}^{\imath\sigma-(n-1)/2-2}.
$$
One expects that this operator does depend on the choice of the caps
$X_\pm$, and thus it is important to understand this dependence. In
particular, one would ideally like to replace these conditions
depending on the caps by boundary conditions at $\pa X_0$. Now,
$P_\emptyset^{-1}$ is characterized by $P_\emptyset^{-1}\psi$, $\psi\in\CI(X)$, having only the $(\mu+\imath
0)^{\imath\sigma}$-type conormal behavior at $\mu=0$, not the
$(\mu-\imath 0)^{\imath\sigma}$ behavior, namely having the form
$$
(\mu+\imath 0)^{\imath\sigma}b_++b_-,
$$
with $b_\pm$ smooth,
since $(\mu+\imath
0)^{\imath\sigma}$ has wave front set in the sink, where the dual
variable $\xi$ of $\mu$ is positive. Restricting to $X_0$ near
the joint boundary $Y_+$ with $X_+$, this has the form
$$
x_0^{2\imath\sigma} a^+_{X_0,+}+a^-_{X_0,+},
$$
with $a^\pm_{X_0,+}$
smooth (and even), while restricting to $X_+$ near $Y_+$ we get the
form
$$
x_+^{2\imath\sigma} a^+_{X_+}+a^-_{X_+},
$$
where, with tilde denoting restriction to $Y_+$,
$$
\tilde a^-_{X_0,+}=\tilde b_-=\tilde a^-_{X_+},\ \tilde
a^+_{X_0,+}=e^{-\pi\sigma}\tilde b_+=e^{-\pi\sigma}\tilde a^+_{X_+}.
$$
Thus, if $\phi$ is supported in $X_0$,
$x_{X_0}^{-\imath\sigma+(n-1)/2}P_\emptyset^{-1}x_{X_0}^{\imath\sigma-(n-1)/2-2}\phi$
is a generalized eigenfunction of $\Delta_{X_+}-(n-1)^2/4-\sigma^2$
with asymptotic behavior
$$
x_+^{\imath\sigma+(n-1)/2} a^+_{X_+}+x_+^{-\imath\sigma+(n-1)/2}a^-_{X_+},
$$
with the result that
$$
\tilde a^-_{X_+}=\cS_{X_+}(\sigma)\tilde a^+_{X_+},
$$
where $\cS_{X_+}(\sigma)$ is the scattering matrix of the
asymptotically hyperbolic problem. In terms of $X_0$ we thus have
$$
\tilde a^-_{X_0,+}=e^{\pi\sigma}\cS_{X_+}(\sigma)\tilde a^+_{X_0,+}.
$$
Since a similar statement also holds at $Y_-$, this Feynman propagator
corresponds to the non-local boundary conditions
$$
a^-_{X_0,\pm}|_{Y_\pm}=e^{\pi\sigma}\cS_{X_\pm}(\sigma)a^+_{X_0,\pm}|_{Y_\pm},
$$
where all $\pm$ signs are consistent on this line. The anti-Feynman
propagator on $X$ produces $(\mu-\imath 0)^{\imath\sigma}$ type conormal
distributions, with the result that
$$
a^-_{X_0,\pm}|_{Y_\pm}=e^{-\pi\sigma}\cS_{X_\pm}(\sigma)a^+_{X_0,\pm}|_{Y_\pm},
$$
then. It would then be an interesting question to study these boundary
conditions directly, as well as more general boundary conditions where
the scattering matrices are replaced by more general
pseudodifferential operators on $Y_\pm$ of order $-2\imath\sigma$,
perhaps even simply $\Delta_{Y_\pm}^{-\imath\sigma}$, which would give
a canonical propagator even in this case. Of course, if one wants to
use a pseudodifferential operator that {\em is} actually the
scattering matrix for a suitable asymptotically hyperbolic space, one
is set!

\section*{Acknowledgments}
I would like to thank Jan Derezi\'nski, Christian G\'erard, Richard
Melrose, Valter Moretti, Michal
Wrochna and Maciej Zworski for helpful discussions and Jesse
Gell-Redman for comments on an earlier version of the manuscript. In particular the
subject of this paper was brought to my attention by
Christian G\'erard.

\bibliographystyle{plain}
\bibliography{sm}

\end{document}